\documentclass[12pt]{article}
\usepackage{amsmath}
\begin{document}
\title{On magic squares\footnote{Delivered to the St. Petersburg Academy
October 17, 1776. Originally published as
\emph{De quadratis magicis}, Commentationes arithmeticae \textbf{2} (1849),
593-602, and reprinted in \emph{Leonhard Euler, Opera Omnia}, Series 1:
Opera mathematica, Volume 7, Birkh\"auser, 1992.
A copy of the original text is available
electronically at the Euler Archive, at www.eulerarchive.org. This paper
is E795 in the Enestr\"om index.}}
\author{Leonhard Euler\footnote{Date of translation: December 4, 2004.
Translated from the Latin by Jordan
Bell, 2nd year undergraduate in Honours Mathematics, School of Mathematics
and Statistics, Carleton University, Ottawa, Ontario, Canada. Email:
jbell3@connect.carleton.ca.
Part of this translation was written
during an NSERC USRA supervised by Dr. B. Stevens.}}
\date{}

\maketitle

1. It is customary for a square to be called a magic square when its
cells are inscribed with the natural numbers in such a way that the
sums of the numbers through each row, column and both diagonals are mutually
equal. Then, if the square were divided into $x$ equal
parts, there would be $xx$ cells altogether, and each of the
rows, columns and both diagonals would contain $x$ cells, in which
each one of the natural numbers $1,2,3,4,\ldots,xx$ would be
arranged, such that the sums for all these lines would
be equal to each other. For this, the sum of all the numbers
from 1 to $xx$ is
\[ \frac{xx(1+xx)}{2}, \]
and the sum for each line is equal to
\[ \frac{x(1+xx)}{2}, \]
by which, for $x=3$ the sum for a single line would be equal to 15.

2. Thus into however many cells a square is divided, the sum of
the numbers deposited in each line can be easily calculated, from which
the sums for all the lines of each square can themselves be easily
calculated.

\begin{tabular}{p{3cm}|p{3cm}|p{3cm}}
$x$&$xx$&$\frac{x(1+xx)}{2}$\\
\hline
1&1&1\\
2&4&5\\
3&9&15\\
4&16&34\\
5&25&65\\
6&36&111\\
7&49&175\\
8&64&260\\
9&81&369\\
&etc.&
\end{tabular}

where $x$ denotes the number of equal parts across into which the square is divided,
and $xx$ the number of cells contained in the square, and where $\frac{1}{2}x(1+xx)$
indicates the sum of all the numbers contained in each line.

3. To help us investigate a certain rule for forming magic squares
of all orders, it is very interesting to observe that all the numbers
1, 2, 3 to $xx$ can be represented with this formula:
\[ mx+n, \]
for if we successively have $m$ take all the values
0, 1, 2, 3, 4 to $x-1$, and then $n$ take all the values
$1,2,3,4,\ldots,x$, it is clear that all the numbers from 1 to $xx$
can be represented by combining each of the values of $m$ with
each of the values of $n$. Furthermore, all the numbers to be inscribed on the
square with the formula $mx+n$ are able to be expressed using two parts,
always in this order, where we use the Latin letters $a,b,c,d$ etc.
for the first part $mx$, and the Greek letters $\alpha,\beta,\gamma,
\delta$ etc. for the second part $n$, where it is clear that for
any number $x$, there is always a Latin and Greek letter that can be
equal to $x$ by having the Latin letters be $0x,1x,2x,3x$ to
$(x-1)x$ and the Greek letters take the values $1,2,3,4,\ldots,x$.
However, this ordering of the Latin and Greek letters
is not fixed, and any Latin letter can denote $0x,1x,2x$ etc., as
long as all the different values are taken by them, with the same
holding for the Greek letters.

4. Now, any number that is to be inscribed on the square can be
represented with a pair of a Latin and a Greek letter,
say by $b+\delta$ or $a+\beta$, etc.; that is, each number
can be represented with these two parts. If each of the Latin letters
are joined together with each of the Greek letters, it is clear
that all the numbers from 1 to $xx$ should result, and it is also
clear that each different combination of letters always produces
different numbers, with no number able to be repeated.

5. Therefore all the numbers are able to be represented by the combination
of a Latin and Greek letter. This in fact yields a rule for the construction of
magic squares. First, the Latin letters are inscribed
in every cell of the square so that the sum for all the lines
are the same, where if there were $x$ letters there would be $xx$
cells altogether in the square, and it is clear that each letter would
be repeated $x$ times. Similarly, the Greek letters are then inscribed
in all the cells of the square, such that the sums for all the lines
are equal. Then, for all the lines the sums of these numbers
made by a combination of a Latin and Greek letter
would be equal. Furthermore, in an arrangement where every Latin letter
is combined with every Greek letter, with this method none of the numbers from from
1 to $xx$ would be missed, and neither would any of them occur
twice.

6. For using this rule to make each type of square according to
how many cells it has, it is clear right away to start with nine
cells, because a square with four cells does not have enough room
for such an arrangement. Furthermore, in general it is seen that
for each type there are $x$ Latin and Greek letters, and that
all the lines have the same number of cells, with the
given conditions satisfied if each line has all the Latin and Greek
letters inscribed in it. However, if the same letter occurs two
or three times in some line, it is always then necessary to have
the sum of all the letters occurring in each line equal to the sums
of all the Latin letters $a+b+c+d+$ etc. or Greek letters
$\alpha+\beta+\gamma+\delta+$ etc.

\begin{center}
{\Large I. Types of squares divided into 9 cells}
\end{center}

7. For this type, it follows that $x=3$, and we have the Latin letters
$a,b,c$ and the Greek letters $\alpha,\beta,\gamma$, where
the Latin letters have the values 0, 3, 6 and the Greek letters
1, 2, 3. We now begin with the Latin letters
$a,b,c$, and it is easy to inscribe them in the 9 cells of our
square such that in each row and column all three letters occur. For
instance, this figure can be seen:
\begin{displaymath}
\begin{matrix}
a&b&c\\
b&c&a\\
c&a&b
\end{matrix}
\end{displaymath}
where now too in one diagonal each of the three letters $a,b,c$
appears, but in the other the same letter $c$ is repeated three
times; it is easy to see that it is not possible to have all the three
different letters in both of the lines at once. However, this does not
cause a problem as long as the diagonal $3c$ is equal to the
sum of the other diagonal $a+b+c$; that is, providing that
$2c=a+b$. From this, it is clear that $c$ should be taken to be
3, and the letters $a$ and $b$ assigned the values 0 and 6; thus
it would be $2c=a+b$. Hence it would be possible to have either
$a=0$ or $b=0$, and from this, the sum of each line results
as $a+b+c=9$.

8.  Similarly, the Greek letters can be distributed in a square
of this type, and we can represent them in this figure with
an inverse arrangement:
\begin{displaymath}
\begin{matrix}
\gamma&\beta&\alpha\\
\alpha&\gamma&\beta\\
\beta&\alpha&\gamma
\end{matrix}
\end{displaymath}
for which it is necessary to have $2\gamma=\alpha+\beta$, and thus
$\gamma=2$. Then, if we combine in a natural way each of the cells from the first
figure with each of the cells from the second figure, it will
be seen that each of the Latin letters is combined with each of the
Greek letters, such that from this combination all the numbers
from 1 to 9 result; this would produce the following figure:
\begin{displaymath}
\begin{matrix}
a\gamma&b\beta&c\alpha\\
b\alpha&c\gamma&a\beta\\
c\beta&a\alpha&b\gamma
\end{matrix}
\end{displaymath}
where it is noted that two letters being joined together does not
mean the product, but instead denotes a combination.

9. With it taken in this figure that $c=3$ and $\gamma=2$,
then the letters $a$ and $b$ can be assigned 0 and 6, and also
the letters $\alpha$ and $\beta$ the values 1 and 3. If we suppose
that $a=0$ and $b=6$, and that $\alpha=1$ and $\beta=3$, the following
magic square will be seen:
\begin{displaymath}
\begin{matrix}
2&9&4\\
7&5&3\\
6&1&8
\end{matrix}
\end{displaymath}
where the sum for each line is 15. If we permute the values
of the letters $a$ and $b$, and likewise $\alpha$ and $\beta$,
it is clear that a different square will be seen. 

10. It is clear that this is a sufficient arrangement of Latin
and Greek letters, and of particular importance in this
is the placement such that each Latin letter is combined with each
Greek letter, and in our arrangement this seems to have
occurred by coincidence. So that we do not have to leave this
to coincidence, before proceding we observe that the arrangement of the
Greek letters $\alpha,\beta,\gamma$ does not depend on the arrangement
of the Latin letters $a,b,c$. Thus for each line, it could be set
that the Greek letters are combined with their Latin equivalents,
e.g. $\alpha$ with $a$, $\beta$ with $b$ and $\gamma$ with $c$. Hence
the first row could be set $a\alpha,b\beta,c\gamma$, and since
the same Greek letter would not occur twice in any row
or column, it can simply be for the second row $b\gamma, c\alpha,
a\beta$, and for the third $c\beta,a\gamma,b\alpha$, from which
this square results:
\begin{displaymath}
\begin{matrix}
a\alpha&b\beta&c\gamma\\
b\gamma&c\alpha&a\beta\\
c\beta&a\gamma&b\alpha
\end{matrix}
\end{displaymath}
where, because in the left diagonal the Greek letter $\alpha$
occurs three time, it is necessary that $3\alpha=\alpha+\beta+\gamma$,
that is, $2\alpha=\beta+\gamma$, which then determines the
value of $\alpha$, namely $\alpha=2$. In the way we see
that $c=3$. However, this does not make a new magic square.

11. Although in this first type the arrangement of the Greek letters
is not difficult to carry out, for squares with larger numbers of cells
it is useful to give a method by which to inscribe
the Greek letters after the Latin letters have been deposited.
For this, a line is chosen in the middle of the rows, columns or
diagonals, such that on either side of the line, the cells that
are equally far away contain different Latin letters. For instance,
in this case the middle column is taken, around which in
the first row are the letters $a$ and $c$, in the second $b$ and $a$,
and in the third $c$ and $b$, such that two different letters
are always on each side.

12. Then after such a middle line has been chosen, in it
each Latin letter is combined with its Greek equivalent, and then
on either side of this, the Greek letters are placed with their
reflected equivalents; for instance, here in this way such a figure
results:
\begin{displaymath}
\begin{matrix}
a\gamma&b\beta&c\alpha\\
b\alpha&c\gamma&a\beta\\
c\beta&a\alpha&b\gamma
\end{matrix}
\end{displaymath}
in which we have clearly combined all the Latin letters with all
the Greek letters. Then, so that the conditions can be satisfied
for the diagonal, we take that
\[ 2c=a+b \textrm{ and } 2\gamma=\alpha+\beta. \]
This figure is in fact not different from the one which
we made earlier in {\S}8. Also, it can be seen that no matter how
the rows and columns are permuted, the sums for the rows and columns
are not changed. However, for the diagonals this can make a very large
difference; if the first column were taken and put on the other
end, this figure would be seen:
\begin{displaymath}
\begin{matrix}
b\beta&c\alpha&a\gamma\\
c\gamma&a\beta&b\alpha\\
a\alpha&b\gamma&c\beta
\end{matrix}
\end{displaymath}
where for the diagonals it must be taken
\[ 2a=b+c \textrm{ and } 2\beta=\alpha+\gamma, \]
for which it is noted that everything has been transposed, which
is an observation that will be very helpful for the following types.

\begin{center}
{\Large II. Types of squares divided into 16 cells}
\end{center}

13. Here it is $x=4$, and we have the four Latin letters
$a,b,c,d$, with the values 0, 4, 8, 12, and also the four
Greek letters $\alpha,\beta,\gamma,\delta$ with the four values
1, 2, 3, 4. Therefore we first inscribe these Latin letters in the
square, such that in each row and column all the four letters
occur, and if it is possible, to have this in both diagonals also.

14. Since there is no prescribed arrangement for these letters
$a,b,c,d$,
in the first row we inscribe them in order, and for the
left diagonal, in the second cell of this diagonal we place
either the letter $c$ or $d$. If we were to write $c$, all the
other letters would then be determined, providing that it is made
sure that the same letter is not written twice in any row
or column. We form the following figure in this way:
\begin{displaymath}
\begin{matrix}
a&b&c&d\\
d&c&b&a\\
b&a&d&c\\
c&d&a&b
\end{matrix}
\end{displaymath}
where each diagonal contains all four letters, so that no
conditions are prescribed for the values of the letters $a,b,c,d$.
Also, if in place of the second cell of this diagonal
we were to write the letter $d$, the resulting figure would be the
only other different possible arrangement, such that with this figure all the other possibilities
would be taken care of.

15. Now, for inscribing the Greek letters, since no middle row or
column is given, we take the diagonal $a,c,d,b$ as the center, and we
find at once that in the cells equally far away on either side
the two letters are distinct, from which it is seen that the rule given before
in \S 11 can be used. Therefore, first we combine the letters in this diagonal
with their Greek equivalents, and then combine the Greek letters with
their reflected equivalents; in this way, the following figure is formed:
\begin{displaymath}
\begin{matrix}
a\alpha&b\delta&c\beta&d\gamma\\
d\beta&c\gamma&b\alpha&a\delta\\
b\gamma&a\beta&d\delta&c\alpha\\
c\delta&d\alpha&a\gamma&b\beta
\end{matrix}
\end{displaymath}

16. Thus in this figure, all the four Latin and Greek letters occur in
all the rows, columns and the full diagonals, and because of this the four
values of these letters can be set without any restrictions. Since there
are 24 variations of four letters, altogether 576 different figures can be
formed, and a good many of the ones made in this way have
structures that are mutually different.

17. By no means should it be thought that all types of magic squares can
be made according to this figure. There are many others that can be made,
where each row does not contain all four Latin and Greek letters,
and that nevertheless fulfill the prescribed conditions. 
Some of these can be made by transposing columns or rows; for instance,
if in the above figure the first column is put at the end, this figure
will be seen:
\begin{displaymath}
\begin{matrix}
b\delta&c\beta&d\gamma&a\alpha\\
c\gamma&b\alpha&a\delta&d\beta\\
a\beta&d\delta&c\alpha&b\gamma\\
d\alpha&a\gamma&b\beta&c\delta
\end{matrix}
\end{displaymath}
where indeed in each row and column all the
Latin and Greek letters still appear, but where in the left
diagonal, descending to the right only two Latin letters occur,
namely $b$ and $c$, and in which also is only a pair of
Greek letters, $\alpha$ and $\delta$. In the other diagonal
are the other two Latin letters $a$ and $d$, and as before,
the Greek letters $\alpha$ and $\delta$.

18. So that this figure satisfies the prescribed conditions,
each letter can take no more than a single value, which
suggests the restriction for the Latin letters of:
\[ b+c=a+d, \]
and similarly for the Greek letters that:
\[ \alpha+\delta=\beta+\gamma; \]
and so, if we take $a=0$, then it follows at once $d=12$,
and that $b=4$ and $c=8$, or vice versa, that $c=4$ and $b=8$.
Similarly for the Greek letters, if we take $\alpha=1$, then it
must be $\delta=4$, and then $\beta=2$ and $\gamma=3$. From this,
a magic square determined as such is made:
\begin{displaymath}
\begin{matrix}
8&10&15&1\\
11&5&4&14\\
2&16&9&7\\
13&3&6&12
\end{matrix}
\end{displaymath}
where clearly the sum for each line is 34. 
Indeed, a great many other forms are able to be made
in this way, by transposing rows or columns.

19. Neither is it absolutely required that each row or column
have all the Latin and Greek letters occuring in it, as rows
and columns can be made with only two Latin and Greek letters,
providing that the sum of them is the same as the sum of all
four letters. It is indeed useful to construct figures with this
special method, and with much work a particular rule can
be made for placing each Latin letter with each Greek letter, such
that while there is not a single sum for all the lines, each Latin
letter is still combined with each Greek letter.

20. To give an example of this method, first we set that
\[ a+d=b+c, \]
and we place the Latin letters as follows,
\begin{displaymath}
\begin{matrix}
a&a&d&d\\
d&d&a&a\\
b&b&c&c\\
c&c&d&d
\end{matrix}
\end{displaymath}
where clearly for all the lines the sum of the numbers is the same.
Then, the Greek letters are combined with each of their Latin
equivalents in the left diagonal, since the two corresponding
letters placed
on either side of this line are different, and then the
Greek letters are combined with their reflected equivalents,
according to which the following figure is made: 
\begin{displaymath}
\begin{matrix}
a\alpha&a\delta&d\beta&d\gamma\\
d\alpha&d\delta&a\beta&a\gamma\\
b\delta&b\alpha&c\gamma&c\beta\\
c\delta&c\alpha&b\gamma&b\beta
\end{matrix}
\end{displaymath}
where for the Greek letters it is necessary that it be taken:
\[ \alpha+\delta=\beta+\gamma; \]
thus, if we take $a=0,b=4,c=8,d=12$ and $\alpha=1,\beta=2,\gamma=3$
and $\delta=4$, such a magic square is seen:
\begin{displaymath}
\begin{matrix}
1&4&14&15\\
13&16&2&3\\
8&5&11&10\\
12&9&7&6
\end{matrix}
\end{displaymath}

21. In this way, many other figures can be formed, for instance
the following:
\begin{displaymath}
\begin{matrix}
a\alpha&d\beta&a\delta&d\gamma\\
b\delta&c\gamma&b\alpha&c\beta\\
d\alpha&a\beta&d\delta&a\gamma\\
c\delta&b\gamma&c\alpha&b\beta
\end{matrix}
\end{displaymath}
where it is clear for the Latin letters to take
\[ a+d=b+c, \]
and for the Greek letters,
\[ \alpha+\delta=\beta+\gamma, \]
from which, if we took the values from above, the following magic
square will be seen:
\begin{displaymath}
\begin{matrix}
1&14&4&15\\
8&11&5&10\\
13&2&16&3\\
12&7&9&6
\end{matrix}
\end{displaymath}

22. In all of these forms, the sums of the Latin and Greek letters
come to the same sum. Others can also be made, which do not follow
any pattern, where the same sum
for all the values is still obtained, but it would be futile
to consider these anomalies, because chance plays such a great
part with them that no fixed pattern can be given for them, and
so
in the following types, chance will be particularly kept in mind,
so that the values of the Latin and Greek letters are not
restricted. 

\begin{center}
{\Large III. Types of squares divided into 25 cells}
\end{center}

23. Therefore, for this type the five Latin letters
$a,b,c,d,e$ occur, and the five Greek letters
$\alpha,\beta,\gamma,\delta,\varepsilon$,  
for which the values of the former are 0, 5, 10, 15, 20, and
the values of the latter 1, 2, 3, 4, 5; both of these letters must
be inscribed in the cells of the square in an arrangement such that
all the letters occur in each row, column and both diagonals.

24. First, therefore, we inscribe all the Latin letters in order
in the top row of the square, and then we fill up the left diagonal
with letters such that the same letter does not occur twice in any of the
remaining lines, with there being more than one way for this to be done.
Once this line has been made, the other diagonal is immediately
determined, and this next figure can be seen:
\begin{displaymath}
\begin{matrix}
a\varepsilon&b\delta&c\gamma&d\beta&e\alpha\\
e\beta&c\alpha&d\delta&a\gamma&b\varepsilon\\
d\alpha&e\gamma&b\beta&c\varepsilon&a\delta\\
b\gamma&d\varepsilon&a\alpha&e\delta&c\beta\\
c\delta&a\beta&e\varepsilon&b\alpha&d\gamma
\end{matrix}
\end{displaymath}
Then below the middle cell is written $a$ and above it $d$, from
which the middle column is completed, and then the remaining lines
are determined immediately.

25. For the Greek letters it is not helpful to use one of the diagonals,
but if we consider the middle column, we find that there are different
letters in the corresponding cells on both sides of it, so thus
we write the Greek equivalents of the Latin letters in this column,
and the Greek letters in the places of their reflected Latin
equivalents, which is how we made this figure.

26. Clearly no restrictions are prescribed for this figure,
and in fact the Latin and Greek letters can take any number.
Since with five letters there are 120 possible permutations,
here altogether 14400 variations can be seen.

27. If we permute the rows or columns between themselves, we can
obtain many other forms, for which we must however set
certain dependencies on the diagonals; for instance, if the first
column is placed at the end, the following form is seen:
\begin{displaymath}
\begin{matrix}
b\delta&c\gamma&d\beta&e\alpha&a\varepsilon\\
c\alpha&d\delta&a\gamma&b\varepsilon&e\beta\\
e\gamma&b\beta&c\varepsilon&a\delta&d\alpha\\
d\varepsilon&a\alpha&e\delta&c\beta&b\gamma\\
a\beta&e\varepsilon&b\alpha&d\gamma&c\delta
\end{matrix}
\end{displaymath}
where in fact in all the rows and columns all the letters occur;
for the diagonals to be satisfied at the same time, this sum
\[ 3c+b+d+3\delta+\beta+\varepsilon \]
and this one
\[ 3a+b+c+3\varepsilon+\alpha+\beta \]
are set to be equal to the sum of all the Latin and Greek letters, namely
\[ a+b+c+d+e+\alpha+\beta+\gamma+\delta+\varepsilon, \]
and when collected, these two equations follow:
\[ 2c+2\delta=a+e+\alpha+\gamma \]
and
\[ 2a+2\varepsilon=d+e+\gamma+\delta ,\]
whose conditions could be satisfied in many ways, where indeed the
the Latin and Greek letters can be determined such that
\[ 1) \quad 2c=a+e, \quad 2) \quad 2a=d+e, \quad 3) \quad 2\delta=\alpha+\gamma
\quad \textrm{and} \quad 4) \quad 2\varepsilon=\gamma+\delta. \]
It is clear that the first two of these are satisfied if the 
letters $d,b,a,c,e$ constitute an arithmetic progression, where
it would be obtained that
\[ d=0,b=5,a=10,c=15 \textrm{ and } e=20; \]
and the two remaining conditions are met if the Greek letters
in the order $\alpha,\beta,\delta,\varepsilon,\gamma$ procede
in arithmetic progression, with it obtained that
\[ \alpha=1,\beta=2,\delta=3,\varepsilon=4 \textrm{ and }
\gamma=5, \]
from which such a square will be seen:
\begin{displaymath}
\begin{matrix}
8&20&2&21&14\\
16&3&15&9&22\\
25&7&19&13&1\\
4&11&23&17&10\\
12&24&6&5&18
\end{matrix}
\end{displaymath}

where clearly all the sums are equal to 65.

28. However, distributing the letters is
by no means simple work, and requires careful consideration.
In particular for the above types where many elements remain
at our discretion, the number of such figures is very large;
removing a restriction causes much work, because there is no
clear restriction prescribed for the values of the letters. If
the letter $c$ takes the middle value, which is 10, with the others
remaining at our discretion, we can fill one diagonal with the 
letter $c$, from which the other letters follow naturally,
and such a figure can be seen:
\begin{displaymath}
\begin{matrix}
c&d&e&a&b\\
b&c&d&e&a\\
a&b&c&d&e\\
e&a&b&c&d\\
d&e&a&b&c
\end{matrix}
\end{displaymath}
Now in the middle row, each of the Greek letters are written
with their Latin equivalents, and then on either side of this
the Greek equivalents are reflected,
such that
the following form is seen:
\begin{displaymath}
\begin{matrix}
c\delta&d\varepsilon&e\alpha&a\beta&b\gamma\\
b\varepsilon&c\alpha&d\beta&e\gamma&a\delta\\
a\alpha&b\beta&c\gamma&d\delta&e\varepsilon\\
e\beta&a\gamma&b\delta&c\varepsilon&d\alpha\\
d\gamma&e\delta&a\varepsilon&b\alpha&c\beta
\end{matrix}
\end{displaymath}
from which $\gamma$ clearly takes the middle value, which is 3;
if we choose the following ordering:
\[ a=0, b=5, c=10, d=15, c=20 \]
and
\[ \alpha=1,\beta=2,\gamma=3,\delta=4,\varepsilon=5, \]
the following magic square will be seen:
\begin{displaymath}
\begin{matrix}
14&20&21&2&8\\
10&11&17&23&4\\
1&7&13&19&25\\
22&3&9&15&16\\
18&24&5&6&12
\end{matrix}
\end{displaymath}

29. By the given method for forming odd squares by switching the values,
this figure is formed:
\begin{displaymath}
\begin{matrix}
11&24&7&20&3\\
4&12&25&8&16\\
17&5&13&21&9\\
10&18&1&14&22\\
23&6&19&2&15
\end{matrix}
\end{displaymath}
which, with it restricted by our formulas, we first consider
the left diagonal in for $c=10$, and take
\[ \delta=1,\alpha=2,\gamma=3,\varepsilon=4 \textrm{ and }
\beta=5, \]
and then
\[ b=0,d=20,a=15,e=5, \]
where from these values this square is made.

30. It is possible to discover other types than the forms
that satisfy these rules, and it is indeed possible to thus
greatly increase the number of magic squares. Yet it is
hardly ever possible to be certain that we have exhausted
all the possibilities, although the number of them is certainly
not infinite. Certainly without doubt, such an investigation for
finding a more general rule for using in different situations
would still not work in many cases. However, it would still
be very beautiful to add to the theory of combinations such a method.

\begin{center}
{\Large IV. Types of squares divided into 36 cells}
\end{center}

31. Since the number of variations here is exceedingly large
and there are many determinations remaining at our discretion,
we produce a specific rule here, with which the Latin and Greek
letters can easily be arranged, where the six Latin letters
take such values:
\[ a+f=b+e=c+d \]
and similarly for the Greek letters,
\[ \alpha+\zeta=\beta+\varepsilon=\gamma+\delta; \]
and by its similarity to \S 20, in each row, we inscribe two Latin
letters, and then in each of the columns in the same way, we place 
two Greek letters, and in this way the following figure is
obtained:\footnote{Translator: As noted in  the reprint of this
paper
in the \emph{Opera Omnia} by the editor, this arrangement
does not indeed make a magic square, as in one diagonal $b\beta$ is
placed twice, and in the other diagonal $e\varepsilon$ appears
twice. However,
the editor
notes that in
Euler's {\em Recerches sur une nouvelle esp\`ece de quarr\'es magiques},
Verhandelingen uitgegeven door het zeeuwsch Genootschap
der Wetenschappen te Vlissingen \textbf{9} (1782), 85-239, reprinted
in the same volume of the \emph{Opera Omnia} as this paper, Euler does in
fact give the
following magic square with 36 cells:
\[ \begin{matrix}3&36&30&4&11&27\\22&13&35&12&14&15\\16&18&8&31&17&21\\
28&20&6&29&19&9\\32&23&25&2&24&5\\10&1&7&33&26&34\end{matrix} \]
}
\begin{displaymath}
\begin{matrix}
a\alpha&a\zeta&a\beta&f\varepsilon&f\gamma&f\delta\\
f\alpha&f\zeta&f\beta&a\varepsilon&a\gamma&a\delta\\
b\alpha&b\zeta&b\beta&e\varepsilon&e\gamma&e\delta\\
e\zeta&e\alpha&e\varepsilon&b\beta&b\delta&b\gamma\\
c\zeta&c\alpha&c\varepsilon&d\beta&d\delta&d\gamma\\
d\zeta&d\alpha&d\varepsilon&c\beta&c\delta&c\gamma
\end{matrix}
\end{displaymath}

32. Thus it can be easily seen that the letters can be arranged
in all the even types successfully,
  and for the odd types this can
be done with the method described earlier, in which the letter
that takes the middle value is repeated in one of the diagonals
and in the other the letters are arranged appropriately. Therefore,
for however many cells a given square has, it is always in our
power to construct many magic squares, even if these rules
that have been given are particular.

\end{document}